\documentclass{amsart}
\usepackage[linesnumbered]{algorithm2e}
\usepackage{amsfonts}
\usepackage{amsmath}
\usepackage{graphicx}
\usepackage{pdfpages}
\usepackage{xcolor}
\usepackage[justification=centering]{caption}
\usepackage[backend=bibtex, sorting=none, bibstyle=alphabetic, citestyle=alphabetic, sorting=nyt,maxbibnames=99]{biblatex}
\bibliography{references.bib}
\AtNextBibliography{\footnotesize}

\newcommand{\Li}{\textrm{li}}
\newcommand{\LLi}{\textrm{Li}}
\renewcommand{\Re}{\textrm{Re}}

\newtheorem{theorem}{Theorem}
\newtheorem{proposition}[theorem]{Proposition}

\numberwithin{theorem}{section}
\theoremstyle{remark}

\begin{document}
\title[Legendre and arithmetic bias]{The Legendre approximation and arithmetic bias}
\author[G. Hiary \and M. Paasche]{Ghaith Hiary \and Megan Paasche}

\address{
    MP: Middlebury College, 
    14 Old Chapel Rd
    Box 3772,
    Middlebury, VT 05753, USA
}
\email{mpaasche@middlebury.edu}
\address{
    GH: Department of Mathematics, The Ohio State University, 231 West 18th
    Ave, Columbus, OH 43210, USA
}
\email{hiary.1@osu.edu}
\subjclass[2010]{Primary: 11-03, 01A55.}
\keywords{Prime number theorem, Legendre, arithmetic bias}
\thanks{* We thank Filip Paun and Jessica Tsang, whose generous donation partly funded this research. We thank the Undergraduate Research program ROMUS in the Mathematics Department at the Ohio State University, Columbus, for their support.}
\maketitle

\begin{abstract}
An interesting episode in the history of the prime number theorem concerns a formula proposed by Legendre for counting the primes below a given bound. We point out that arithmetic bias likely played an important role in arriving at that formula and in its subsequent widespread, decades-long recognition. We also show that the Legendre constant $1.08366$ satisfies a certain simple and natural criterion, and conjecture that this criterion is how Legendre arrived at that erroneous constant in his formula.
\end{abstract}

\section{Introduction}

Let $\pi(x)$ denote the number of primes $\le x$. A central result in number theory is the prime number theorem, which gives an asymptotic formula for the prime-counting function $\pi(x)$. Some of Euler's unpublished notes, found by Bouniakowsky and Chebyshev and included in \textit{Commentationes Arithmeticae},
show that Euler considered the problem of counting the primes \cite{euler_1969}. He tabulated the primes up to $x=2,4,8,\ldots,1024$ and observed that the ``law of the primes'' appears to be that they become less frequent as $x$ increases. Euler further remarked, however, that the number of primes in centuries (blocks of length $100$) does not decrease regularly, but changes in a very irregular manner. He gave isolated examples around higher $x$ (near $90{,}000$) to demonstrate this.\footnote{Euler might have had an intuitive understanding of the prime number theorem. Pertinent evidence for this can be found in a letter he wrote to Goldbach in 1752 \cite{euler_goldbach_1752}, 
though it is hard to discern Euler's precise intention there. Considering Euler's famous work on the convergent sum $\zeta(2)=\sum_n n^{-2}$, in comparison with the harmonic sum over primes $\sum_p 1/p$, which Euler proved diverges, Euler at least knew that the primes occur more frequently than the squares of integers, so $\pi(x)\ge \sqrt{x}$ for $x$ large enough. Moreover, Euler surely knew the primes are sparse among all integers considering that he proved, using his mathematical language at the time, and frequently cited the result that $\sum_p 1/p$ diverges far more slowly than the usual harmonic series.}

Euler's observation about the irregular distribution of the primes manifests as difficulty in approximating $\pi(x)$ accurately using a simple continuous function. Nevertheless, in 1808, Legendre published such an approximation of $\pi(x)$ in his book \textit{Essai Sur la Théorie des Nombres}.
His attempt was based on a more extensive tabulation of the primes than Euler had, going up to $x=400{,}000$ at first (2nd edition) \cite{legendre_1808} and extended to $x=1{,}000{,}000$ by 1830 (3rd edition) \cite{legendre_1830}.  
He conjectured that $\pi(x)$ is well-approximated by the function
\begin{equation}\label{legendre_formula}
\frac{x}{A\log x-B}, 
\end{equation}
where $A=1$ and $B=1.08366$. But he also acknowledged that the distribution of the primes is subject to anomalies.  Moreover, in his concluding remarks, Legendre expressed astonishment that basic analysis, combined with his approximation of $\pi(x)$, enabled accurate evaluation of sums and products over the primes. This remark may be seen as a precursor to the subject of analytic number theory.

Several prominent mathematicians at the time cast doubt on the significance of the constant $1.08366$ in the Legendre formula-- specifically, Gauss \cite{gauss_1849}, Dirichlet \cite{dirichlet_1838}, and Chebyshev \cite{chebyshev_1852}. 
They favored the approximation of $\pi(x)$ given by logarithmic integral, 
\begin{equation}\label{li_asymp}
\begin{split}
\Li(x) &= \int_2^x \frac{1}{\log t}\,dt \\
&\sim \frac{x}{\log x}+\frac{1!\, x}{(\log x)^2}+\cdots+\frac{\ell!\, x}{(\log x)^{\ell+1}}+\cdots, 
\end{split}
\end{equation}
where the asymptotic in the second line is as $x\to\infty$. See \cite{davenport_1967}.\footnote{The closely related integral $\LLi(x)$ has the same definition as $\Li(x)$ except one starts integrating at $x=0$ and uses the Cauchy principal value to define the integral around $x=1$. The functions $\Li(x)$ and $\LLi(x)$ differ only by the constant $\LLi(2)=1.04516\ldots$.}

Of course, we now know that $\Li(x)$ is the better approximation asymptotically. Even as early as 1852, Chebyshev proved that 
\begin{equation}\label{alt_limit}
   \frac{x}{\pi(x)}-\log x 
\end{equation}
cannot have a limit different from $-1$ as $x\to \infty$
\cite{chebyshev_1852}, which is consistent with the second term in the asymptotic expansion \eqref{li_asymp}.
In comparison, the Legendre formula gives a limit equal to $-1.08366$. Therefore, if $\pi(x)$ is indeed well-approximated by the expression \eqref{legendre_formula} with $A=1$, then the Legendre constant $1.08366$ cannot be the best choice for $B$ asymptotically. For a short proof of the Chebyshev result concerning the limiting value of \eqref{alt_limit}, see \cite{pintz_1980}.

Other than the leading term $x/\log x$, the Legendre formula disagrees with the more natural-looking approximation $\Li(x)$. It is remarkable that Legendre published his choice of $B$ extended to 5 decimal places, indicating a high degree of confidence in his conjecture.\footnote{A confidence interval of $\pm 0.00001$ for the value of $B$ presumably reflects an expectation that the approximation \eqref{legendre_formula} with $A=1$ is typically reliable to within $\pm 1$ around $x=1{,}000{,}000$.} It is also surprising that the odd-looking number $1.08366$ did not appear connected to any known mathematical constants, despite the fundamental nature of $\pi(x)$. Drach \cite{drach_1844}
attempted to remedy this apparent imperfection, proposing in 1844 that the correct constant should be $\log(5\sqrt{\pi}/3) = 1.08319\ldots$ This became moot just a few years later, though, after Chebyshev published his work. 

In Proposition~\ref{B_prop}, and the discussion following it, we show that the Legendre constant $1.08366$ satisfies a certain simple and natural criterion. We conjecture that this criterion is how Legendre arrived at that constant.

\section{Arithmetic bias}

After careful examination of Legendre's original work, as well as of contemporary work from that period, we could not find a direct explanation for why Legendre chose the approximation model \eqref{legendre_formula} and the values $A=1$ and $B=1.08366$, beyond that these choices were based on numerical tables in some way. We provide a likely explanation for the origin of Legendre's choices for \eqref{legendre_formula}.
Along the way, we show that Legendre was likely misled by a conspiracy of arithmetic bias that steered him towards the approximation \eqref{legendre_formula}, and he otherwise made a reasonable choice of $B$ considering the range of primes tables available to him. 

The arithmetic bias in question is the inequality
\begin{equation}\label{bias}
    \Li(x)>\pi(x),
\end{equation}
which typically holds, and always does so in the range $8\le x\le 1{,}000{,}000$. See \cite{shanks_1959} and \cite{rubinstein_sarnak_1994} for discussions of various types of arithmetic bias. 

The bias \eqref{bias} is not as artificial as it may seem at first, and certainly goes counter to basic probabilistic considerations.
It is unclear whether Gauss considered the inequality \eqref{bias} to be highly probable, but Riemann 
\cite{riemann_1859}
indicated belief that $\LLi(x)$ ought to give a value of $\pi(x)$ that is slightly too large. This is not always the case, however, as  Littlewood~\cite{littlewood_1914} established in 1914 that $\Li(x)-\pi(x)$ changes sign infinitely often.

\begin{figure}[h!]
    \includegraphics[scale = 0.48]{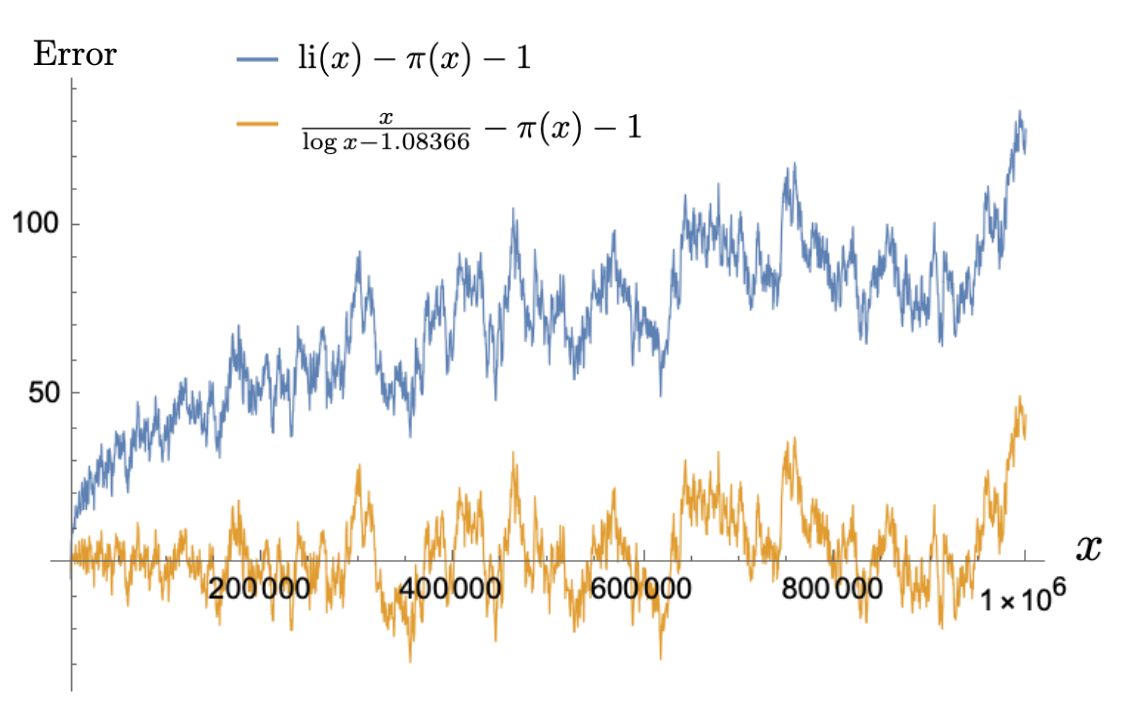}
    \caption{\small Illustration of arithmetic bias and the Legendre approximation for $x\le 10^6$.}
    \label{arithm_bias_1}
\end{figure}

Theoretical justification for the arithmetic bias \eqref{bias} stems from Riemann's 1859 paper and his exact analytic formulas for $\pi(x)$ and related functions.\footnote{In Riemann's paper, $\pi(x)$ denotes the number of primes $< x$ if $x$ is not a prime, and if $x$ is a prime then $\pi(x)$ is equal to its average value from the left and right. This definition differs slightly (it is smaller by $1/2$) from the definition in the introduction when $x$ is exactly equal to a prime.} After applying the Mobius inversion, Riemann obtained
\begin{equation}\label{pi_formula}
    \pi(x)=\sum_{n=1}^{\infty} \frac{\mu(n)}{n}f(x^{1/n}),
\end{equation}
where $\mu$ is the Mobius function and $f(x)$ is defined as in Riemann's 1859 paper. So,
\begin{equation}
    f(x)=\sideset{}{'}\sum_{n\le x} \frac{\Lambda(n)}{\log n},
\end{equation}
where $\Lambda$ is the von Mangoldt function and the primed sum means that if $x$ is a prime power then $f(x)$ is defined by its average from the left and the right.

Riemann gave an exact analytic formula for $f(x)$. He expresses $f(x)$ as a main term equal to $\LLi(x)$, followed by an infinite sum of secondary terms over the nontrivial zeros $\rho$ of Riemann's zeta $\zeta(s)$, plus other small and easily understood terms. Furthermore, the secondary terms corresponding to the $\rho$'s may be taken to be of the form $\LLi(x^{\rho})+\LLi(x^{1-\rho})$. See \cite[p. 48]{edwards_1974} for details. 

The $n=2$ summand on the right-side of \eqref{pi_formula} corresponds to the contribution of the squares of primes, and is about $-\sqrt{x}/\log x$. Thus, the contribution of the squares of primes is negative and noticeable. Although this contribution could be overwhelmed by the contributions of the $\rho$'s from the $n=1$ summand, this usually does not happen. This is because the largest contributor among the $\rho$'s is the first zero $\rho_1 = 1/2+14.134725\ldots i$, and its contribution is limited in size by about $2\sqrt{x}/(|\rho_1|\log x)$. Hence, since $\rho_1$ has a high ordinate (imaginary part), the contribution of $\rho_1$ will be significantly less than that of the squares of primes. 

While there are infinitely more contributions from the $\rho$'s to account for, these contributions, assuming the Riemann hypothesis that $\Re(\rho)=1/2$ for all $\rho$, are of decreasing size. Furthermore, they are oscillatory, with oscillations that are expected to be independent. Therefore, their sum can be expected to exhibit substantial cancellation. Overall, then, the $-1$ from the squares of primes should usually win over the $2/|\rho_1|$, and this materializes as the arithmetic bias \eqref{bias}.\footnote{Therefore, the arithmetic bias should disappear if instead of $\pi(x)$ and $\Li(x)$ one considers $f(x)$ and its approximation by $\LLi(x)$ -- that is, if one counts the primes with suitable weights.}

It seems highly unlikely that
Legendre would have been aware of the arithmetic bias \eqref{bias} when he made his conjecture.
Once Legendre decided that the approximation for $\pi(x)$ should be of the form \eqref{legendre_formula}, it remained to choose the constants $A$ and $B$. The choice $A=1$ was probably predetermined, in a sense. 
Gauss had already guessed by 1793 that the density of primes around $x$ is about $1/\log(x)$. Although this  was unknown to Legendre at the time\footnote{And indeed,  Legendre's work on $\pi(x)$ was apparently unknown or forgotten by Gauss, as indicated in Gauss' 1849 letter to his former student Encke \cite{gauss_1849}.}, it appears to have reached the status of a folklore conjecture by the beginning of the 19th century. 
For example, Dirichlet \cite{dirichlet_1838}, while expressing skepticism about the accuracy of the Legendre approximation, wrote in a footnote in 1838 that the true approximation of $\pi(x)$ is $\sum_{2\le n\le x}1/\log n \sim \Li(x)$, which in view of the asymptotic expansion \eqref{li_asymp} dictates that $A=1$.

With $A$ chosen to be $1$, one could choose $B$ so as to minimize the average error in the approximation \eqref{legendre_formula}. For a given $x$, this error is
\begin{equation}
    \mathcal{E}(x,B) = \frac{x}{\log x - B} - \pi(x)-1, 
\end{equation}
where the extra $-1$ is because Legendre counted $1$ as a prime. Alternatively, one could average the difference in \eqref{alt_limit} over $x$ and choose $B$ accordingly. We will examine the first method, as it makes the connection with arithmetic bias more transparent, and it returns a value of $B$ closer to what Legendre obtained.

\begin{figure}[h!]
    \includegraphics[scale = 0.42]{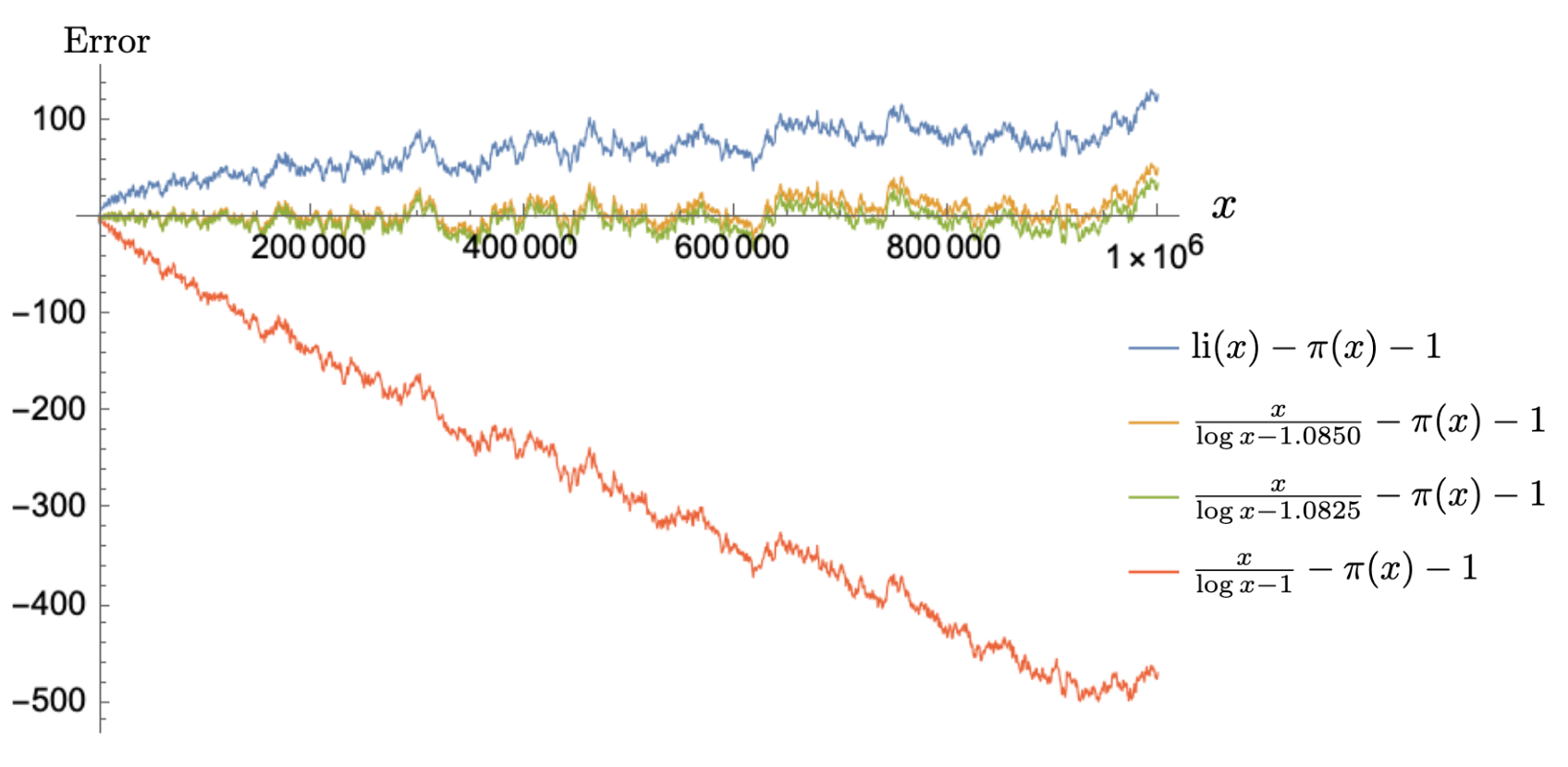}
    \caption{\small Tracking the error for various $B$ and for $\Li(x)$ over $x\le 10^6$.}
    \label{error_plot}
\end{figure}

Let $\langle \mathcal{E}(x,B) \rangle_x$ be the average of $\mathcal{E}(x,B)$ over the primes $x$ with $3\le x\le 10^6$.\footnote{We start at $x=3$ to ensure that the denominator $\log x - B$ in the Legendre approximation is positive. When $x=2$, this denominator is negative if $B > \log 2 = 0.693…$, as we will have.}
This is the same range of primes that was available to Legendre.
In Table~\ref{avg_error},  we give the value of $\langle \mathcal{E}(x,B) \rangle_x$ for various $B$, showing that the
value of $B$ that minimizes $|\langle \mathcal{E}(x,B) \rangle_x|$ is between $1.0825$ and $1.0850$. 
Moreover, Figure~\ref{error_plot} tracks $\mathcal{E}(x,B)$ for $B=1, 1.0825, 1.0850$ as $x$ ranges over $3\le x\le 10^6$, as well as tracking the error produced by the approximation $\Li(x)$ over the same range of $x$. In view of Figure~\ref{error_plot}, the approximation $\Li(x)$ certainly appears as an unlikely candidate to approximate $\pi(x)$. The approximation by $\Li(x)$ produces an error that is heavily biased, large, and positive. In comparison, the approximation \eqref{legendre_formula} with $A=1$ and $B\in[1.0825,1.0850]$ appears to be the far better approximation. These approximations give errors that behave more intuitively, moderately fluctuating about $0$.

\begin{table}[ht]
\renewcommand\arraystretch{1.5}
\centering
\begin{tabular}{|c|c||c|c|}
\hline
 $B$ & $\langle \mathcal{E}(x,B) \rangle_x $ & $B$ & $\langle \mathcal{E}(x,B) \rangle_x$\\
 \hline
 1.0700 & -41.2565  & 1.0825 & -0.95052 \\
 1.0725 & -33.202 & 1.0850 & 7.12087    \\
 1.0750 & -25.1442  & 1.0875 & 15.1958  \\
 1.0775 & -17.083 & 1.0900 & 23.2744    \\
 1.0800 & -9.01846 & 1.0925 & 31.3572   \\
\hline
\end{tabular}
\caption{\small The average error for various $B$.}
\label{avg_error}
\end{table}

Since for all $x$ the expression \eqref{legendre_formula} is monotonically increasing in $B <\log x$, the average error $\langle \mathcal{E}(x,B) \rangle_x$ is increasing in $B < \log 3 = 1.0986\ldots$ (since we start the averaging at $x=3$). So, in view of Table~\ref{avg_error} which shows that $\langle \mathcal{E}(x,B) \rangle_x$ assumes both positive and negative values, we deduce that there is a unique $B$ that makes $|\langle \mathcal{E}(x,B) \rangle_x|$ equal to zero. Therefore, we obtain the following proposition.

\begin{proposition}\label{B_prop}
There is a unique real number $B_0 <\log 3$ that makes $\langle \mathcal{E}(x,B_0) \rangle_x$ vanish. The number $B_0$ satisfies  $1.0825 \le B_0 \le 1.0850$.
\end{proposition}

One can further  refine the value of $B$ using a bisection method. The bisection method allows us to numerically search for this unique root $B$. To this end, we computed $\langle \mathcal{E}(x,B) \rangle_x$ for several values of $B$ in an interval around $1.0825$ and recursively chose a smaller interval around the value of $B$ that reduces $\left|\langle \mathcal{E}(x,B) \rangle_x\right|$. This process results in choosing 
$$B=1.08279,$$ 
which is close to Legendre’s constant $1.08366$. The difference is equal to merely $0.00087$. 

There may be a few reasons why Legendre chose a slightly different constant than we found. For example, he indicates in \textit{Essai Sur la Théorie des Nombres} that he might have used his own unpublished table of primes. Interestingly, Legendre mentions that he calculated a value for $\pi(1{,}000{,}000)$ that differs from the value in the W\'ega, Chernac, and Burckhardt tables cited in his book. 
And both values differ from the true value: by $-28$ using Legendre's calculation, and by $6$ using the cited tables. Such discrepancies easily explain the $0.00087$ difference that we found. 

We therefore speculate that Legendre used a nearly equivalent method to the minimization method we described in this section to choose $B$. This is evidenced by how well we are able to capture the Legendre constant using this minimization procedure. Thus, by design, the Legendre approximation yields much better results than $\Li(x)$ and better results than approximations with other values of $B$. The superior accuracy of the Legendre approximation in the range of $x$ he considered is probably the reason his approximation remained prominent for several decades, despite the inexplicable constant $1.08366$. From an empirical viewpoint, the arithmetic bias \eqref{bias} was making it hard to recognize $\Li(x)$ as the true approximation in the feasible range of $x$ at that time.



%

\section{Limitations of numerical data}

Any reasonable investigator, faced with a heavily biased error, may abandon the logarithmic integral approximation $\Li(x)$ unless they had special knowledge about the arithmetic bias phenomenon \eqref{bias}. Legendre favored the approximation \eqref{legendre_formula} with $A=1$ and then calibrated $B$, thus resolving the bias issue in the range he considered. In contrast, the approximation $\Li(x)$ did not offer any parameters to calibrate. It took Riemann's 1859 paper to realize that all the calibration in that approximation occurs via the positions of the nontrivial zeta zeros $\rho$.

But even in the absence of arithmetic bias, so assuming $\LLi(x)$ is a very accurate approximation for $\pi(x)$, if one insists on approximating $\pi(x)$ using \eqref{legendre_formula} with $A=1$, then one would be led further astray from the asymptotically correct choice $B=1$. Under this hypothetical, the choice of $B$ that minimizes the size of the average error over the primes $x$ with $5\le x\le 10^6$ is even larger, around $1.10407$. 
The other attractive choices of $B$, such as $B=1$ or $B=0$, also result in heavy bias, though in the opposite direction, over the range Legendre considered.   

\begin{figure}[ht]
    \includegraphics[scale = 0.5]{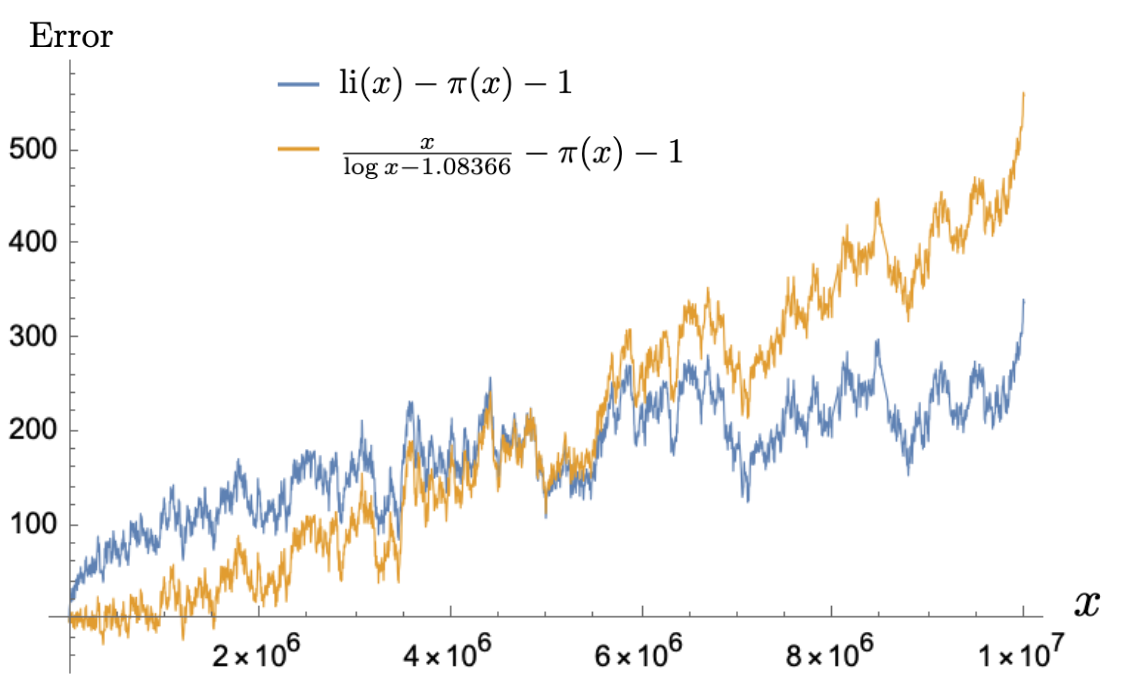}
    \caption{\small Arithmetic bias for $x\le 10^7$.}
    \label{arithm_bias_2}
\end{figure}

So, in retrospect, the data available to Legendre was simply too limited to discern that the true approximation was $\Li(x)$ or that $B=1$ was the asymptotically correct choice. Although the Legendre approximation is more accurate than $\Li(x)$ up to around $x=3{,}000{,}000$, it starts to lose advantage after that, becoming noticeably
worse once $x$ exceeds $6{,}000{,}000$. This is illustrated in Figure~\ref{arithm_bias_2}. 

Therefore, access to more data, like Gauss indicates he had, would have easily cast doubt on the approximation \eqref{legendre_formula}.
Gauss also offered a numerical argument against the Legendre approximation. Gauss \cite{gauss_1849} pointed out that the error in the Legendre approximation, while smaller than that of $\Li(x)$ for $x\le 10^6$, grows much faster. (This is seen in Figure~\ref{arithm_bias_2} as the orange curve is steeper than the blue curve.) So it is conceivable that the error will get too large eventually. Indeed, the error in the Legendre approximation ultimately exceeds that of $\Li(x)$ by a large margin.

\section{Conclusions}
The Legendre approximation and the Legendre constant tell a story about the advantages and limitations of experimental mathematics, and its many twists, turns, and complexities even in most fundamental settings.
Legendre put substantial weight on the numerics and the guidance that resulted from the numerics, and he subsequently gave an asymptotically erroneous value of $B$ to several decimal places. 
While he was likely misled by arithmetic bias and limited data, the Legendre formula is remarkably accurate over the range he considered. This is probably the reason that the Legendre formula became so prevalent in the literature, possibly overshadowing $\Li(x)$ which suffered from arithmetic bias. The Legendre formula was supported by seemingly convincing, though in retrospect limited, numerical evidence, and stood prominent for several decades despite the fair criticisms of Legendre's contemporaries and despite the strange-looking constant $1.08366$. Moreover, Legendre appears to have used the simple and natural criterion of minimizing the average error to choose $B$, thus ensuring the accuracy of his approximation for $x\le 10^6$ by design. 

\newpage 
\printbibliography

\end{document}